\documentclass{amsart}

\newtheorem{theorem}{Theorem}[section]
\newtheorem{lemma}[theorem]{Lemma}

\theoremstyle{definition}

\usepackage{amscd,amssymb}
\begin{document}

\title[Defining Relations of Invariants of two $3 \times 3$ matrices]
{Defining Relations of Invariants \\
of two $3 \times 3$ matrices}

\author[Helmer Aslaksen, Vesselin Drensky and Liliya Sadikova]
{Helmer Aslaksen, Vesselin Drensky and Liliya Sadikova}
\address{Department of Mathematics, National University of Singapore,
          Singapore 117534, Singapore}
\email{aslaksen@math.nus.edu.sg}
\address{Institute of Mathematics and Informatics,
          Bulgarian Academy of Sciences,
          1113 Sofia, Bulgaria}
\email{drensky@math.bas.bg}
\address{Fachgruppe Informatik,
          RWTH Aachen, 52056 Aachen, Germany}
\email{sadikova@stce.rwth-aachen.de}

\thanks{The second author was partially supported by Grant MM1106/2001
of the Bulgarian National Science Fund.}
\subjclass{16R30}

\begin{abstract}
Over a field of characteristic 0, the algebra of invariants of several
$n\times n$ matrices under simultaneous conjugation by $GL_n$
is generated by traces of products of generic matrices. Teranishi, 1986,
found a minimal system of eleven generators of the algebra of invariants of two
$3\times 3$ matrices. Nakamoto, 2002, obtained an explicit, but very complicated,
defining relation for a similar system of generators over $\mathbb Z$.
In this paper we have found another natural set of eleven generators of this algebra
of invariants over a field of characteristic 0 and have given the defining relation
with respect to this set. Our defining relation is much simpler than that
of Nakamoto. The proof is based on easy computer calculations
with standard functions of Maple but the explicit form of the relation has been found
with methods of representation theory of general linear groups.
\end{abstract}

\maketitle

\section*{Introduction}

Let $K$ be any field of characteristic 0 and let
$X_i=\left(x_{pq}^{(i)}\right)$, $p,q=1,\ldots,n$, $i=1,\ldots,d$, be $d$
generic $n\times n$ matrices. The conjugation of $X_i$
with the invertible $n\times n$ matrix $g$,
\[
X_i=\left(x_{pq}^{(i)}\right)\to gX_ig^{-1}=\left(y_{pq}^{(i)}\right),
\]
defines an action of the general linear group $GL_n=GL_n(K)$ on the
polynomial algebra in $n^2d$ variables
\[
\Omega_{nd}=K[x_{pq}^{(i)}\mid p,q=1,\ldots,n,\ i=1,\ldots,d]
\]
by $g\ast x_{pq}^{(i)}=y_{pq}^{(i)}$, $g\in GL_n$. We denote by $C_{nd}$
the algebra of invariants $\Omega_{nd}^{GL_n}$.
Its elements are called the invariants under the action of $GL_n$ by
simultaneous conjugation of $d$ matrices of size $n\times n$.

Traditionally, we call a result giving the explicit generators of the
algebra of invariants of a linear group $G$ a
first fundamental theorem of the invariant theory of $G$
and a result describing the relations between the generators a
second fundamental theorem.

For a background on the algebra of matrix invariants see e.g., \cite{F2, DF},
for the history see \cite{F3}.
It is known that $C_{nd}$ is generated by traces of products of generic matrices
$\text{\rm tr}(X_{i_1}\cdots X_{i_k})$. For a fixed $n$ an upper bound for
the degree $k$ is given in terms of PI-algebras. By the Nagata-Higman theorem
the nil algebras of bounded index are nilpotent, i.e.,
the polynomial identity $x^n=0$ implies the identity $x_1\cdots x_m=0$.
Then $k\leq m$ and for $d$ sufficiently large this bound is sharp.
A description of the defining relations of $C_{nd}$ is given in the theory of
Razmyslov \cite{R} and Procesi \cite{P} in the language of ideals of the group
algebras of symmetric groups.

Explicit minimal sets of generators of $C_{nd}$ are known for $n=2$ and any $d$
and for $n=3$ and $d=2$.
By a theorem of Sibirskii \cite{S}, $C_{2d}$ is generated
by $\text{\rm tr}(X_i)$, $1\leq i\leq d$, $\text{\rm tr}(X_iX_j)$, $1\leq i\leq j\leq d$,
$\text{\rm tr}(X_iX_jX_k)$, $1\leq i<j<k\leq d$. There are no relations between
the five generators of the algebra $C_{22}$, i.e.,
$C_{22}\cong K[z_1,\ldots,z_5]$. For $d=3$, Sibirskii \cite{S} found one relation
and Formanek \cite{F1} proved that all relations follow from it.
Aslaksen, Tan and Zhu \cite{ATZ} found, mostly by computer,
some defining relations for bigger $d$,
as well as generators and some defining relations
for the $2\times 2$ matrix invariants for other classical groups.
The centre of $GL_2$ acts trivially on $\Omega_{2d}$ and we have also a
natural action of $PSL_2$ on $\Omega_{2d}$. Since $PSL_2({\mathbb C})$ is isomorphic to
$SO_3({\mathbb C})$, we may apply invariant theory of orthogonal groups, see
Procesi \cite{P} and Le Bruyn \cite{L}. In particular,
Drensky \cite{D2} translated the defining relations of the invariants of
$SO_3$ and obtained the defining relations of $C_{2d}$ for all $d$.

Teranishi \cite{T} found the following system of generators of $C_{32}$:
\begin{equation} \label{Teran-gen}
\begin{array}{c}
\text{\rm tr}(X),\text{\rm tr}(Y),\text{\rm tr}(X^2),\text{\rm tr}(XY),\text{\rm tr}(Y^2), \\
\\
\text{\rm tr}(X^3),\text{\rm tr}(X^2Y),\text{\rm tr}(XY^2),\text{\rm tr}(Y^3),
\text{\rm tr}(X^2Y^2),\text{\rm tr}(X^2Y^2XY),
\end{array}
\end{equation}
where $X$, $Y$ are generic $3 \times 3$ matrices.
He showed that the first ten of these generators
form a homogeneous system of parameters of $C_{32}$ and
$C_{32}$ is a free module with generators 1 and $\text{\rm tr}(X^2Y^2XY)$
over the polynomial algebra on these ten elements.
(Teranishi \cite{T} found also a set of generators and
a homogeneous system of parameters of $C_{42}$.)
Abeasis and Pittaluga \cite{AP} found a system of generators of
$C_{3d}$ in terms of representation theory of the symmetric and general linear groups,
in the spirit of its use in theory of PI-algebras.

The algebras $C_{nd}$ have a natural multigrading which takes into account the degrees
of the traces $\text{\rm tr}(X_{i_1}\cdots X_{i_k})$ with respect to each of the
generic matrices $X_1,\ldots,X_d$. The Hilbert series of $C_{nd}$
is defined as the formal power series
\[
H(C_{nd},t_1,\ldots,t_d)=\sum_{k_i\geq 0}
\dim(C_{nd}^{(k_1,\ldots,k_d)})t_1^{k_1}\cdots t_d^{k_d}
\]
with coefficients equal to the dimensions of the homogeneous components
$C_{nd}^{(k_1,\ldots,k_d)}$ of degree $(k_1,\ldots,k_d)$.
It carries a lot of information for the algebra.
The Hilbert series of $C_{2d}$ are calculated e.g., in \cite{L}
and those of $C_{32}$ and $C_{42}$ in \cite{T}. Van den Bergh \cite{V}
reduced the determination of the Hilbert series of $C_{nd}$ to a
problem about flows in a certain graph and obtained an
important consequence for the denominators of the rational
functions in the explicit form of the series. He also evaluated
$H(C_{nd},t_1,\ldots,t_d)$ for several small $n$ and $d$.

We shall need the Hilbert series of $C_{32}$ obtained by Teranishi \cite{T}
as a consequence of his description of $C_{32}$:
\begin{equation} \label{Teran-series}
H(C_{32},t_1,t_2)=
\frac{1+t_1^3t_2^3}
{(1-t_1)(1-t_2)q_2(t_1,t_2)q_3(t_1,t_2)(1-t_1^2t_2^2)},
\end{equation}
where the commuting variables $t_1$ and $t_2$ count, respectively,
the degrees of $X$ and $Y$ and
\[
q_2(t_1,t_2)=(1-t_1^2)(1-t_1t_2)(1-t_2^2),
\]
\[
q_3(t_1,t_2)=(1-t_1^3)(1-t_1^2t_2)(1-t_1t_2^2)(1-t_2^3).
\]
It follows from the description of the generators of $C_{32}$, that
$\text{\rm tr}(X^2Y^2XY)$ satisfies a quadratic equation with coefficients
depending on the other ten generators. The explicit (but very complicated)
form of the equation was found by Nakamoto \cite{N}, over $\mathbb Z$,
with respect to a slightly different system of generators.
Revoy \cite{Re} studied the field of rational invariants of two $3\times 3$ matrices
and also found eleven generators, with an explicitly given relation between them.

In this paper we have found another natural set of eleven generators of the algebra
$C_{32}$ and have given the defining relation with respect to this set. It has turned out
that our relation is much simpler than that in \cite{N}.

The first change in the set of generators is well known.
In the traces of products of degree $\geq 2$ we have replaced the generic matrices
$X,Y$ with generic traceless matrices $x,y$, so our generators become
$\text{\rm tr}(X),\text{\rm tr}(Y)$ and
$\text{\rm tr}(x^2),\text{\rm tr}(xy),\text{\rm tr}(y^2)$,
$\text{\rm tr}(x^3),\text{\rm tr}(x^2y),\text{\rm tr}(xy^2),\text{\rm tr}(y^3)$,
$\text{\rm tr}(x^2y^2),\text{\rm tr}(x^2y^2xy)$. Extending diagonally
the natural action of $GL_2$ on $K\cdot X+K\cdot Y$ to $C_{32}$, we have used
the result of Abeasis and Pittaluga \cite{AP} and we have replaced
the generators $\text{\rm tr}(x^2y^2),\text{\rm tr}(x^2y^2xy)$
with elements which generate one-dimensional $GL_2$-modules.
It becomes clear that the only defining relation between our generators
spans a one-dimensional $GL_2$-submodule of the subalgebra of $C_{32}$
generated by all traces of products of $x$ and $y$. Using several rules for
computing with tensor products and symmetric algebras of irreducible $GL_2$-modules,
we have seen that the possible candidate for the relation is a linear combination
with unknown coefficients of nine explicitly found elements in $C_{32}$.
In order to determine these coefficients, we have evaluated the hypothetic
relation for various concrete traceless matrices $x,y$. The evaluations are equal to 0, and
we have obtained a linear system with respect to the nine unknown coefficients of the considered linear
combination. Using standard procedures of Maple,
we have obtained the only solution of the system. Then, having the hypothetic
defining relation, it is easy to check that it is really a relation.

The paper is organized as follows. In Section 1 we give our system of generators
and the defining relation between them. Our proof is direct and
does not use any representation theory of $GL_2$. In Section 2 we show how have we
found our generators and the relation. We believe that our methods can be successfully
used for other similar problems.

The main results of the present paper have been announced in \cite{ADS}.

\section{Main Results}
\label{sec-main}

It is a standard trick in the study of matrix invariants to replace the generic matrices
in the traces of products with generic traceless matrices. We express $X$ and $Y$ in
the form
\begin{equation}\label{from ordinary to traceless matrices}
X=\frac{1}{3}\text{\rm tr}(X)e+x\quad \text{\rm and}\quad
 Y=\frac{1}{3}\text{\rm tr}(Y)e+y,
\end{equation}
where $e$ is the identity matrix and $x$, $y$ are generic traceless matrices.
Changing the variables $x_{ii}$ and $y_{ii}$, we may assume that
\begin{equation}\label{matrices}
x=\left( \begin{array}{ccc}
                x_{11} & x_{12} & x_{13} \\
                x_{21} & x_{22} & x_{23} \\
                x_{31} & x_{32} & -(x_{11}+x_{22})\\
             \end{array} \right), \
y=\left( \begin{array}{ccc}
                y_{11} & y_{12} & y_{13} \\
                y_{21} & y_{22} & y_{23} \\
                y_{31} & y_{32} & -(y_{11}+y_{22})\\
             \end{array} \right).
\end{equation}
By well known arguments, as for ``ordinary'' generic matrices,
without loss of generality we may assume that $x$ is a diagonal matrix, i.e.,
\begin{equation}\label{x is diagonal}
x=\left( \begin{array}{ccc}
                x_1 & 0 & 0 \\
                0 & x_2 & 0 \\
                0 & 0 & -(x_1+x_2)\\
             \end{array} \right).
\end{equation}
Till the end of the paper we fix the notation $x,y$ for the two generic traceless matrices.
Replacing $X,Y$ with their expressions from
(\ref{from ordinary to traceless matrices}), we obtain
that, instead of the generators (\ref{Teran-gen}) of Teranishi, the algebra $C_{32}$ is
generated by the system
\begin{equation}\label{new-gen}
\begin{array}{c}
\text{\rm tr}(X),\text{\rm tr}(Y),\text{\rm tr}(x^2),\text{\rm tr}(xy),\text{\rm tr}(y^2), \\
\\
\text{\rm tr}(x^3),\text{\rm tr}(x^2y),\text{\rm tr}(xy^2),\text{\rm tr}(y^3),
\text{\rm tr}(x^2y^2),\text{\rm tr}(x^2y^2xy).\\
\end{array}
\end{equation}
Now we replace the traces $\text{\rm tr}(x^2y^2),\text{\rm tr}(x^2y^2xy)$
with elements more convenient for our purposes.
We define
\begin{equation}\label{v}
v=\text{\rm tr}(x^2y^2)-\text{\rm tr}(xyxy),
\end{equation}
\begin{equation}\label{w}
w=\text{\rm tr}(x^2y^2xy)-\text{\rm tr}(y^2x^2yx).
\end{equation}

The following lemma justifies the introduction of $v$ and $w$.
It expresses $\text{\rm tr}(x^2y^2)$ and $\text{\rm tr}(x^2y^2xy)$ in terms of
$v,w$ and the other generators.

\begin{lemma} \label{vw-lemma}
The following equations hold in $C_{32}$:
 \begin{equation}\label{trace of xxyy}
   \text{\rm tr}(x^2y^2)=\frac{1}{3}v+\frac{1}{6}\text{\rm tr}(x^2)\text{\rm tr}(y^2)
                         + \frac{1}{3}\text{\rm tr}^2(xy),
 \end{equation}

\begin{equation}\label{trace of xxyyxy}
\begin{array}{c}
   \text{\rm tr}(x^2y^2xy)=\frac{1}{2}w+\frac{1}{6}\text{\rm tr}(xy)v
        +\frac{1}{12}\text{\rm tr}(x^2)\text{\rm tr}(xy)\text{\rm tr}(y^2) \\
\\
        +\frac{1}{6}\text{\rm tr}^3(xy)-\frac{1}{6}\text{\rm tr}(x^3)\text{\rm tr}(y^3)
            +\frac{1}{2}\text{\rm tr}(x^2y)\text{\rm tr}(xy^2). \\
\end{array}
\end{equation}
\end{lemma}

\begin{proof} For the proof of the first equality we shall make use of the
Cayley-Hamilton theorem. For any $3\times 3$ matrix $z$ with eigenvalues
$z_1,z_2,z_3$
\[
  z^3-e_1z^2+e_2z-e_3e=0,
\]
where $e_i$ is the $i$-th elementary symmetric function in $z_1,z_2,z_3$.
The Newton formulas
\[
p_k-e_1p_{k-1}+\cdots+(-1)^{k-1}e_{k-1}p_1+(-1)^kke_k=0,\quad k=1,2,3,
\]
allow to express the $e_i$'s in terms of
the power sums $p_k=z_1^k+z_2^k+z_3^k$:
\[
e_1=p_1,\quad e_2=\frac{1}{2}(p_1^2-p_2),\quad
e_3=\frac{1}{6}(2p_3-3p_1p_2+p_1^3).
\]
Since $p_k=\text{\rm tr}(z^k)$, if $\text{\rm tr}(z)=0$, then we obtain
\[
z^3-\frac{1}{2}\text{tr}(z^2)z-\frac{1}{3}\text{tr}(z^3)e=0.
\]
Multiplying by $z$ and taking the trace we derive the trace identity for
$3\times 3$ traceless matrices
\begin{equation}\label{CHtheorem for traceless matrices}
  \text{\rm tr}(z^4)-\frac{1}{2}\text{\rm tr}^2(z^2)=0.
\end{equation}
Since the trace of the matrix $x+y$ is equal to 0, the identity
(\ref{CHtheorem for traceless matrices}) gives
\[
\text{\rm tr}((x+y)^4)-\frac{1}{2}\text{\rm tr}^2((x+y)^2)=0.
\]
The multihomogeneous components of every trace identity are also trace identities.
Considering the component of second degree with respect to $x$,
and taking into account the invariance of the trace under cyclic permutation
of the variables, we obtain
\[
4\text{\rm tr}(x^2y^2)+2\text{\rm tr}(xyxy)
-\text{\rm tr}(x^2)\text{\rm tr}(y^2)-2\text{\rm tr}^2(xy)=0.
\]
Replacing $\text{\rm tr}(xyxy)$ with $\text{\rm tr}(x^2y^2)-v$
from (\ref{v}), we derive (\ref{trace of xxyy}).
In order to prove (\ref{trace of xxyyxy}), we may proceed in a similar way. Instead, using
Maple, we have evaluated (\ref{trace of xxyyxy}) on the matrices $x$
from (\ref{x is diagonal}) and $y$ from (\ref{matrices}) and
have obtained that both sides are equal.
\end{proof}
\par

Now we define the following elements of $C_{32}$:
\begin{equation}\label{u}
u=\left| \begin{array}{cc}
          \text{\rm tr}(x^2)& \text{\rm tr}(xy) \\
          \text{\rm tr}(xy) & \text{\rm tr}(y^2) \\
          \end{array} \right|,
\end{equation}
\begin{equation}\label{w1 and the others}
 w_1=u^3,  \quad
 w_2=u^2v, \quad
 w_4=uv^2, \quad
 w_7=v^3,
\end{equation}
\begin{equation}\label{w5}
  w_5=v\left| \begin{array}{ccc}
                \text{\rm tr}(x^2)& \text{\rm tr}(xy) & \text{\rm tr}(y^2) \\
                \text{\rm tr}(x^3)& \text{\rm tr}(x^2y) & \text{\rm tr}(xy^2) \\
                \text{\rm tr}(x^2y) & \text{\rm tr}(xy^2) & \text{\rm tr}(y^3) \\
             \end{array} \right| ,
\end{equation}
\begin{equation}\label{w6}
 w_6=\left| \begin{array}{cc}
          \text{\rm tr}(x^3)& \text{\rm tr}(xy^2) \\
          \text{\rm tr}(x^2y) & \text{\rm tr}(y^3) \\
          \end{array} \right|^2
   -4\left| \begin{array}{cc}
          \text{\rm tr}(y^3)& \text{\rm tr}(xy^2) \\
          \text{\rm tr}(xy^2) & \text{\rm tr}(x^2y) \\
          \end{array} \right|\left| \begin{array}{cc}
          \text{\rm tr}(x^3)& \text{\rm tr}(x^2y) \\
          \text{\rm tr}(x^2y) & \text{\rm tr}(xy^2) \\
          \end{array} \right|,
\end{equation}
\begin{equation}\label{w31}
w_3'=u\left| \begin{array}{ccc}
                \text{\rm tr}(x^2)& \text{\rm tr}(xy) & \text{\rm tr}(y^2) \\
                \text{\rm tr}(x^3)& \text{\rm tr}(x^2y) & \text{\rm tr}(xy^2) \\
                \text{\rm tr}(x^2y) & \text{\rm tr}(xy^2) & \text{\rm tr}(y^3) \\
             \end{array} \right| ,
\end{equation}
\begin{equation}\label{w32}
\begin{array}{c}
  w_3''=5[\text{\rm tr}^3(y^2)\text{\rm tr}^2(x^3)
+\text{\rm tr}^3(x^2)\text{\rm tr}^2(y^3)]\\
\\
   -30[\text{\rm tr}^2(y^2)\text{\rm tr}(xy)\text{\rm tr}(x^2y)\text{\rm tr}(x^3)+
        \text{\rm tr}^2(x^2)\text{\rm tr}(xy)\text{\rm tr}(y^3)\text{\rm tr}(xy^2)]\\
\\
  +3\{[4\text{\rm tr}(y^2)\text{\rm tr}^2(xy)+\text{\rm tr}^2(y^2)
            \text{\rm tr}(x^2)][3\text{\rm tr}^2(x^2y)
+2\text{\rm tr}(xy^2)\text{\rm tr}(x^3)]\\
\\
   +[4\text{\rm tr}^2(xy)\text{\rm tr}(x^2)+\text{\rm tr}^2(x^2)\text{\rm tr}(y^2)]
   [3\text{\rm tr}^2(xy^2)+2\text{\rm tr}(x^2y)\text{\rm tr}(y^3)]\}\\
\\
  -2[2\text{\rm tr}^3(xy)+3\text{\rm tr}(x^2)\text{\rm tr}(xy)\text{\rm tr}(y^2)]
  [9\text{\rm tr}(xy^2)\text{\rm tr}(x^2y)+\text{\rm tr}(x^3)\text{\rm tr}(y^3)],
\end{array}
\end{equation}
where $v$ is defined in (\ref{v}).
The element $w_3''$ can be expressed in the following
simple way. Recall that a linear mapping $\delta$ of an algebra $R$ is a derivation
if $\delta(rs)=\delta(r)s+r\delta(s)$ for all $r,s\in R$. Every mapping
$\delta_0:\{X,Y\}\to C_{32}$ can be extended to a derivation $\delta$ of $C_{32}$
which commutes with the trace. In particular, defining $\delta(X)=0$, $\delta(Y)=X$,
we obtain a derivation $\delta$ which sends $x$ to 0 and $y$ to $x$. Then
\begin{equation}\label{new definition of w32}
w_3''=\frac{1}{144}\sum_{i=0}^6(-1)^i\delta^i(\text{\rm tr}^3(y^2))
\delta^{6-i}(\text{\rm tr}^2(y^3)).
\end{equation}

The following theorem is the main result of our paper.

\begin{theorem}\label{main theorem}
The algebra of invariants $C_{32}$ of two $3\times 3$ matrices has the following
presentation. It
is generated by
\begin{equation}\label{our generators}
\begin{array}{c}
\text{\rm tr}(X), \text{\rm tr}(Y),
 \text{\rm tr}(x^2), \text{\rm tr}(xy), \text{\rm tr}(y^2),\\
\\
\text{\rm tr}(x^3), \text{\rm tr}(x^2y), \text{\rm tr}(xy^2),
\text{\rm tr}(y^3), v, w\\
\end{array}
\end{equation}
subject to the defining relation
\begin{equation}\label{our defining relation}
  w^2-\left( \frac{1}{27}w_1-\frac{2}{9}w_2+\frac{4}{15}w_3'+\frac{1}{90}w_3''
    +\frac{1}{3}w_4-\frac{2}{3}w_5-\frac{1}{3}w_6-\frac{4}{27}w_7 \right)=0,
\end{equation}
where the elements $v, w, w_1,w_2,w_3',w_3'',w_4,w_5,w_6,w_7$
are given in {\rm (\ref{v}), (\ref{w}), (\ref{w1 and the others}), (\ref{w5}),
(\ref{w6}), (\ref{w31})} and  {\rm (\ref{w32})}.
\end{theorem}

\begin{proof} Replacing the system of generators (\ref{Teran-gen}) with the
system (\ref{new-gen}),
Lemma\ \ref{vw-lemma} insures that the elements (\ref{our generators})
are generators of the algebra $C_{32}$.
The relation (\ref{our defining relation}) is checked with the help of Maple,
where $x$ is as in (\ref{x is diagonal}) and $y$ is as in (\ref{matrices}).

Let us consider the polynomial algebra $P_{11}$
in the 11 variables (\ref{our generators})
written in the form $P_{11}=(S[\text{\rm tr}(X),\text{\rm tr}(Y)])[w]$, where
\begin{equation}\label{the algebra S}
S=K[\text{\rm tr}(x^2), \text{\rm tr}(xy), \text{\rm tr}(y^2),\text{\rm tr}(x^3),
\text{\rm tr}(x^2y), \text{\rm tr}(xy^2), \text{\rm tr}(y^3), v],
\end{equation}
and the principal ideal $I=(f)$ of $P_{11}$ generated by the element from
(\ref{our defining relation})
\[
 f = w^2-\left(\frac{1}{27}w_1-\frac{2}{9}w_2+\frac{4}{15}w_3'+\frac{1}{90}w_3''
    +\frac{1}{3}w_4-\frac{2}{3}w_5-\frac{1}{3}w_6-\frac{4}{27}w_7\right).
\]
The factor algebra $P_{11}/I$ is graded and is a free $S[\text{\rm tr}(X),\text{\rm tr}(Y)]$-module,
freely generated by 1 and $w$.
Therefore, the Hilbert series $H(t_1,t_2)$ of $P_{11}/I$ is
$$
\frac{1+t_1^3t_2^3}
{(1-t_1)(1-t_2)(1-t_1^2)(1-t_1t_2)(1-t_2^2)(1-t_1^3)(1-t_1^2t_2)(1-t_1t_2^2)(1-t_2^3)(1-t_1^2t_2^2)},
$$
and is equal to the Hilbert series
(\ref{Teran-series}) of the algebra $C_{32}$.
Since the algebra $C_{32}$ satisfies the relation $f=0$ from (\ref{our defining relation}),
it is a homomorphic image of $P_{11}/I$ and the coefficients of
its Hilbert series are bounded from above by the coefficients of the Hilbert series
of $P_{11}/I$. Since the Hilbert series of $C_{32}$ and $P_{11}/I$ coincide,
the algebras $C_{32}$ and $P_{11}/I$ are isomorphic and this completes the proof.
\end{proof}

\section{How Did We Find The Relation?}

We refer e.g., to \cite{M} for a background on representation theory of $GL_d$
and to \cite{D1}  for its applications in the spirit of the problems
considered here.
The natural action of the general linear group $GL_d$ on the
vector space with basis $\{x_1,\ldots,x_d\}$
(or $\{X_1,\ldots,X_d\}$) can be extended diagonally on the
free algebra $K\langle x_1,\ldots,x_d\rangle$ and on the algebra generated by
all formal traces
\[
K[\text{\rm tr}(X_{i_1}\cdots X_{i_m})\mid i_j=1,\ldots,d,\ m=1,2\ldots].
\]
This induces also a $GL_d$-action on the algebra $C_{nd}$.
It is well known that $C_{nd}$ is a direct sum of irreducible polynomial $GL_d$-modules.
We denote by $W_d(\lambda)$ the irreducible $GL_d$-module indexed by the partition
$\lambda=(\lambda_1,\ldots,\lambda_d)$. We assume that $W_d(\lambda)=0$ if
$\lambda$ is a partition in more than $d$ parts.

Abeasis and Pittaluga \cite{AP} suggested the following way to describe the generators
of $C_{nd}$. The subalgebra $C_k$ of $C_{nd}$ generated by all traces
$\text{\rm tr}(X_{i_1}\cdots X_{i_m})$ of degree $m\leq k$ is a $GL_d$-submodule of
$C_{nd}$. Let $C_{nd}^{(k+1)}$ be the homogeneous component of degree $k+1$
of $C_{nd}$. Then the intersection $C_k\cap C_{nd}^{(k+1)}$ has a complement
$G_{k+1}$ in $C_{nd}^{(k+1)}$,
which is the $GL_d$-module of the ``new'' generators of degree $k+1$.
Without loss of generality we may assume that $G_{k+1}$ is a submodule of
the $GL_d$-module spanned by all traces $\text{\rm tr}(X_{i_1}\cdots X_{i_{k+1}})$.
Abeasis and Pittaluga performed the calculations in the case of invariants
of $3 \times 3$ matrices and obtained the following results for the $GL_d$-module
$G=G_1\oplus G_2\oplus\cdots$ of the generators of $C_{3d}$
(the partitions in \cite{AP} are given in ``Francophone'' way, i.e., transposed to ours):
\[
G=W_d(1)\oplus W_d(2)\oplus W_d(3)\oplus W_d(1^3)\oplus W_d(2^2)\oplus
W_d(2,1^2)
\]
\[
\oplus W_d(3,1^2)\oplus W_d(2^2,1)\oplus W_d(1^5)\oplus
W_d(3^2)\oplus W_d(3,1^3).
\]
For $d=2$ we may assume that $X_1=X$, $X_2=Y$. Then the algebra $C_{32}$
has a system of generators which spans a $GL_2$-module
\[
W_2(1)\oplus W_2(2)\oplus W_2(3)\oplus W_2(2^2)\oplus W_2(3^2).
\]
The module $W_2(1)$ is spanned by the traces $\text{\rm tr}(X)$ and $\text{\rm tr}(Y)$.
Hence the subalgebra of $C_{32}$ generated by traces of products
of the traceless matrices $x$ and $y$ has a system of generators forming the module
\[
W_2(2)\oplus W_2(3)\oplus W_2(2^2)\oplus W_2(3^2).
\]

Let $U_{kd}$ be the $GL_d$-module spanned by all formal traces
$\text{\rm tr}(X_{i_1}\cdots X_{i_k})$. The structure of $U_{kd}$
is determined by the structure of the $S_k$-module $\text{\rm tr}(V_k)$
of the multilinear elements of $U_{kd}$,
where $S_k$ is the symmetric group of degree $k$ and
$\text{\rm tr}(V_k)$ has a basis
$\{\text{\rm tr}(X_kX_{\sigma(1)}\cdots X_{\sigma(k-1)})
\mid \sigma\in S_{k-1}\}$. The $S_k$-module $\text{\rm tr}(V_k)$ is isomorphic to the
induced module of the trivial one-dimensional module of a cycle of length $k$ in $S_k$.
For our purposes, in the following lemma we give a direct proof of the decomposition of
$U_{k2}$ for $k=2,3,4,6$.

\begin{lemma} \label{lemma-A}
The $GL_2$-modules $U_{k2}$, $k=2,3,4,6$, have the following decompositions:
\[
U_{22}=W_2(2),\quad U_{32}=W_2(3),\quad U_{42}=W_2(4)\oplus W_2(2^2),
\]
\[
U_{62}=W_2(6)\oplus 2W_2(4,2)\oplus W_2(3^2).
\]
\end{lemma}
\begin{proof}
If a polynomial $GL_d$-module $W_d$ has the decomposition
\[
W_d=\bigoplus m(\lambda)W_d(\lambda),\quad m(\lambda)\geq 0,
\]
as a direct sum of irreducible components,
then it has a natural grading and its Hilbert series is
\[
H(W_d,t_1,\ldots,t_d)=\sum m(\lambda)S_{\lambda}(t_1,\ldots,t_d),
\]
where $S_{\lambda}(t_1,\ldots,t_d)$ is the Schur function associated with
$\lambda$. The Hilbert series plays the role of the character of $W_d$
and the multiplicities $m(\lambda)$ can be determined in a unique way
from the Hilbert series. For $d=2$ the Schur functions
have the following simple expression:
\[
S_{(a+b,b)}(t_1,t_2)=(t_1t_2)^b(t_1^a+t_1^{a-1}t_2+\cdots +t_1t_2^{a-1}+t_2^a).
\]
The vector spaces $U_{k2}$, $k=2,3,4,6$, have bases, respectively,
\[
\{\text{\rm tr}(X^2),\text{\rm tr}(XY), \text{\rm tr}(Y^2)\},
\]
\[
\{\text{\rm tr}(X^3), \text{\rm tr}(X^2Y), \text{\rm tr}(XY^2), \text{\rm tr}(Y^3)\},
\]
\[
\{\text{\rm tr}(X^4), \text{\rm tr}(X^3Y), \text{\rm tr}(X^2Y^2),  \text{\rm tr}(XYXY),
\text{\rm tr}(XY^3), \text{\rm tr}(Y^4)\},
\]
\[
\{\text{\rm tr}(X^6), \text{\rm tr}(X^5Y), \text{\rm tr}(X^4Y^2),
\text{\rm tr}(X^3YXY), \text{\rm tr}(X^2YX^2Y),
\]
\[
\text{\rm tr}(X^3Y^3), \text{\rm tr}(X^2Y^2XY),
\text{\rm tr}(X^2YXY^2), \text{\rm tr}(XYXYXY),
\]
\[
\text{\rm tr}(X^2Y^4),
\text{\rm tr}(XY^3XY), \text{\rm tr}(XY^2XY^2), \text{\rm tr}(XY^5), \{\text{\rm tr}(X^6)
\},
\]
Hence the Hilbert series of $U_{k2}$ are
\[
H(U_{22},t_1,t_2)= t_1^2+t_1t_2+t_2^2=S_{(2)}(t_1,t_2),
\]
\[
H(U_{32},t_1,t_2)= t_1^3+t_1^2t_2+t_1t_2^2+t_2^3=S_{(3)}(t_1,t_2),
\]
\[
H(U_{42},t_1,t_2)= t_1^4+t_1^3t_2+2t_1^2t_2^2+t_1t_2^3+t_2^4
\]
\[
=S_{(4)}(t_1,t_2)+S_{(2^2)}(t_1,t_2),
\]
\[
H(U_{62},t_1,t_2)=
t_1^6+t_1^5t_2+3t_1^4t_2^2+4t_1^3t_2^3+3t_1^2t_2^4+t_1t_2^5+t_2^6
\]
\[
=S_{(6)}+2S_{(4,2)}+S_{(3^2)}.
\]
This gives the desired decompositions.
\end{proof}

\begin{lemma}\label{lemma-B}
The $GL_2$-submodules $W_2(2^2)\subset U_{42}$
and $W_2(3^2)\subset U_{62}$ are one-dimensional and are
spanned, respectively, by
\[
V=v(X,Y)=\text{\rm tr}(X^2Y^2)-\text{\rm tr}(XYXY),
\]
\[
W=w(X,Y)=\text{\rm tr}(X^2Y^2XY)-\text{\rm tr}(Y^2X^2YX).
\]
\end{lemma}
\begin{proof}
The Schur function
$S_{(b^2)}(t_1,t_2)=(t_1t_2)^b$ gives that the module $W_2(b^2)$ is one-dimensional.
The module $W_d(\lambda)$ is generated by a
unique, up to a multiplicative constant,
homogeneous element $w_{\lambda}$ of degree $\lambda_i$ with respect to $X_i$,
called the highest weight vector of $W_d(\lambda)$. It is characterized by
the following property, which follows from the results of
De Concini, Eisenbud and Procesi \cite{DEP} and Almkvist, Dicks and Formanek \cite{ADF},
see also Koshlukov \cite{K} for the version which we need. We state it
for two variables and for $U_{k2}$ only.
If $w_{\lambda}(X,Y) \in U_{k2}$ is homogeneous of degree
$(\lambda_1,\lambda_2)$ and $w_{\lambda}(X\vert Y,Z)$ is the partial linearization
of $w_{\lambda}$ in $Y$
(i.e., the homogeneous component of degree 1 in $Z$ of $w_{\lambda}(X,Y+Z)$),
then $w_{\lambda}(X,Y)$
is a highest weight vector for some $W_2(\lambda_1,\lambda_2)$
if and only if $w_{\lambda}(X\vert Y,X)=0$. Equivalently, if $\Delta$ is the derivation
of the algebra generated by formal traces of products of $X$ and $Y$,
and $\Delta$ is defined by $\Delta(X)=0$, $\Delta(Y)=X$, then
$w_{\lambda}(X\vert Y,X)=\Delta(w_{\lambda}(X,Y))=0$. Now the proof of the lemma follows from the
equalities
$\Delta(V)=\Delta(W)=0$
which can be verified by direct calculations. For example,
the partial linearization of $V$ is
\[
v(X\vert Y,Z)=\text{\rm tr}(X^2YZ)+\text{\rm tr}(X^2ZY)-2\text{\rm tr}(XYXZ),
\]
which gives that
\[
\Delta(V)=v(X\vert Y,X)=(1+1-2)\text{\rm tr}(X^3Y)=0.
\]
\end{proof}

In order to find the expressions for $V,W$ one may proceed in the following way,
as shown for $W_2(3^2)$. The generator $W$ of $W_2(3^2)$ is homogeneous of degree
$(3,3)$. Hence it is a linear combination of
$\text{\rm tr}(X^3Y^3)$, $\text{\rm tr}(X^2Y^2XY)$,
$\text{\rm tr}(Y^2X^2YX)$, $\text{\rm tr}(XYXYXY)$
and has the form
\[
    W=\eta_1\text{\rm tr}(X^3Y^3)+\eta_2\text{\rm tr}(X^2Y^2XY)
      +\eta_3\text{\rm tr}(Y^2X^2YX)
    +\eta_4\text{\rm tr}(XYXYXY)
\]
with unknown coefficients $\eta_1,\eta_2,\eta_3,\eta_4$.
The partitial linearizations of the summands are:
\[
\text{\rm for}\quad
\text{\rm tr}(X^3Y^3):\quad \text{\rm tr}(X^3XY^2)+ \text{\rm tr}(X^3YXY) +
    \text{\rm tr}(X^3YYX)
\]
\[
=2\text{\rm tr}(X^4Y^2)+\text{\rm tr}(X^3YXY);
\]
\[
\text{\rm for}\quad
    \text{\rm tr}(X^2Y^2XY):\quad \text{\rm tr}(X^2XYXY)+\text{\rm tr}(X^2YXXY)+
    \text{\rm tr}(X^2YYXX)
\]
\[
    =\text{\rm tr}(X^4Y^2)+\text{\rm tr}(X^3YXY)+\text{\rm tr}(X^2YX^2Y);
\]
\[
\text{\rm for}\quad
    \text{\rm tr}(Y^2X^2YX): \quad \text{\rm tr}(XYX^2YX)+\text{\rm tr}(YXX^2YX)+
    \text{\rm tr}(Y^2X^2XX)
\]
\[
    =\text{\rm tr}(X^4Y^2)+\text{\rm tr}(X^3YXY)+\text{\rm tr}(X^2YX^2Y);
\]
\[
\text{\rm for}\quad
    \text{\rm tr}(XYXYXY): \quad 3\text{\rm tr}(X^3YXY).
\]
The condition $W(X\vert Y,X)=0$ gives
\[
(2\eta_1+\eta_2+\eta_3)\text{\rm tr}(X^4Y^2)
+(\eta_1+\eta_2+\eta_3+3\eta_4)\text{\rm tr}(X^3YXY)
\]
\[
+(\eta_2+\eta_3)\text{\rm tr}(X^2YX^2Y)=0.
\]
Since the coefficients of
$\text{\rm tr}(X^4Y^2),\text{\rm tr}(X^3YXY),\text{\rm tr}(X^2YX^2Y)$
are equal to 0, we obtain the system
\[
   \begin{array}{rll}
     2\eta_1+&\eta_2+\eta_3&=0 \\
     \eta_1+&\eta_2+\eta_3+3\eta_4&=0 \\
     &\eta_2+\eta_3&=0.\\
   \end{array}
\]
Its only solution is $\eta_1=\eta_4=0$, $\eta_3=-\eta_2$. Our element $W$ is obtained
for $\eta_2=1$.

In this way we obtain that the subalgebra of $C_{32}$ generated by products of traceless
matrices is generated by the $GL_2$-modules $W_2(2), W_2(3)$ spanned respectively by
$\{\text{\rm tr}(x^2),\text{\rm tr}(xy),\text{\rm tr}(y^2)\}$ and
$\{\text{\rm tr}(x^3),\text{\rm tr}(x^2y),\text{\rm tr}(xy^2),\text{\rm tr}(y^3)\}$,
and the $GL_2$-modules $W_2(2^2)$ and $W_2(3^2)$ spanned, respectively, by
$v$ and $w$ from (\ref{v}) and (\ref{w}). Translating the results of Teranishi
\cite{T}, the submodules $W_2(2), W_2(3)$ and $W_2(2^2)$ of $C_{32}$ generate the
subalgebra $S$ from (\ref{the algebra S}) and the whole subalgebra generated by
traces of products of traceless matrices is a free $S$-module with basis $\{1,w\}$.
Since $w^2 \in S\oplus wS$, we derive that $w$ satisfies a
relation of the form
\begin{equation}\label{the relation for w}
f=w^2+aw+b=0,
\end{equation}
for some $a,b \in S$.

\begin{lemma}\label{the modules generated by the coefficients}
The $GL_2$-modules generated by the coefficients $a$ and $b$ in
{\rm (\ref{the relation for w})} are isomorphic to $W_2(3^2)$ and
$W_2(6^2)$, respectively.
\end{lemma}

\begin{proof}
The elements $g\in GL_2$ act on $w\in W_2(3^2)$ by
\[
g(w(x,y))=w(g(x),g(y))={\det}^3(g)\cdot w(x,y).
\]
Let us generate a $GL_2$-module by $f=w^2+aw+b$ from (\ref{the relation for w}).
Acting with $g\in GL_2$, we obtain
\[
g(f)={\det}^6(g)\cdot w^2+{\det}^3(g)\cdot g(a)w+g(b).
\]
Hence we have also the relation
\[
g(f)-{\det}^6(g)\cdot f={\det}^3(g)(g(a)-{\det}^3(g)\cdot a)w+(g(b)-{\det}^6(g)\cdot b)=0.
\]
Since 1 and $w$ form a basis of the free $S$-module $S\oplus wS$, we obtain that
$g(a)={\det}^3(g)\cdot a$ and $g(b)={\det}^6(g)\cdot b$ for all $g\in GL_2$.
Hence $a$ and $b$ generate $GL_2$-modules isomorphic, respectively, to
$W_2(3^2)$ and $W_2(6^2)$.
\end{proof}

\begin{lemma}\label{multiplicities for 33 and 66}
In the decomposition of $S$ as a $GL_2$-module
\[
S=\bigoplus m(\lambda_1,\lambda_2)W_2(\lambda_1,\lambda_2)
\]
we have $m(3^2)=0$ and $m(6^2)=8$.
\end{lemma}

\begin{proof}
The multiplicities $m(\lambda_1,\lambda_2)$ in the decomposition of $S$
can be obtained using methods of the recent papers \cite{DG1, DG2}.
(In \cite{DG1} Drensky and Genov calculated the generating function
of the multiplicities for the whole algebra $C_{32}$ and in \cite{DG2} they gave
an easier way to calculate the multiplicities also for other objects similar to $C_{32}$.)
The multiplicities $m(3^2)$ and $m(6^2)$ can also be obtained directly, using
well known rules of representation theory of symmetric and general linear groups.
The first possibility is to find the coefficients of degree 6 and 12 in the Hilbert series
of $S$ and to express them in terms of Schur functions, as in Lemma \ref{lemma-A}.
We shall give one more possibility.
The algebra $S$ is the tensor product of three polynomial algebras, namely
the symmetric algebras
$K[W_2(2)]=K[\text{\rm tr}(x^2),\text{\rm tr}(xy),\text{\rm tr}(y^2)]$,
$K[W_2(3)]=K[\text{\rm tr}(x^3),
\text{\rm tr}(x^2y),\text{\rm tr}(xy^2),\text{\rm tr}(y^3)]$
and $K[W_2(2^2)]=K[v]$. The decomposition
\[
K[W_2(2)]=\bigoplus W_2(2\lambda_1,2\lambda_2)
\]
is well known and follows from the equality
\[
\prod_{1\leq i\leq j\leq d}\frac{1}{1-t_it_j}
=\sum S_{(2\lambda_1,\ldots,2\lambda_d)}(t_1,\ldots,t_d).
\]
The decomposition
\[
K[W_2(2^2)]=\bigoplus_{b\geq 0} W_2(2b^2)
\]
is also well known. Using the multiplicities in the decomposition
of $K[W_2(3)]$ given in \cite{DG1} or calculating the first coefficients
of the Hilbert series
\[
H(K[W_2(3)],t_1,t_2)=\frac{1}{(1-t_1^3)(1-t_1^2t_2)(1-t_1t_2^2)(1-t_2^3)}
=\sum_{m\geq 0}h_{3m}(t_1,t_2),
\]
where $h_k$ is the homogeneous component of degree $k$, or simply counting
the basis elements of degree $\leq 12$ in $K[W_2(3)]$, we obtain
\[
h_3(t_1,t_2)=S_{(3)}(t_1,t_2),
\]
\[
h_6(t_1,t_2)=S_{(6)}(t_1,t_2)+S_{(4,2)}(t_1,t_2),
\]
\[
h_{12}(t_1,t_2)=S_{(12)}(t_1,t_2)+S_{(10,2)}(t_1,t_2)+S_{(9,3)}(t_1,t_2)
+S_{(8,4)}(t_1,t_2)+S_{(6^2)}(t_1,t_2).
\]
This gives the parts of the decomposition of $K[W_2(3)]$ which we need.
Now, the homogeneous component of degree 6 of $S$ is obtained in the following way:
\[
S^{(6)}=(K[W_2(2)])^{(6)}\oplus (K[W_2(2)])^{(2)}\otimes (K[W_2(2^2)])^{(4)}
\oplus (K[W_2(3)])^{(6)}.
\]
Using the Littlewood-Richardson rule, which in the case of $GL_2$ states
\[
W_2(a+b,b)\otimes W_2(c+d,d)\cong \bigoplus_{p=0}^c W_2(a+b+d+p,b+d+c-p),
\quad a\geq c,
\]
we calculate
\[
S^{(6)}=(W_2(6)\oplus W_2(4,2))\oplus
W_2(2)\otimes W_2(2^2)\oplus(W_2(6)\oplus W_2(4,2))
\]
\[
=2W_2(6)\oplus 3W_2(4,2),
\]
and this means that $m(3^2)=0$.
Similarly, we apply the Littlewood-Richardson rule to
\[
S^{(12)}=(K[W_2(2)])^{(12)}\oplus (K[W_2(2)])^{(8)}\otimes (K[W_2(2^2)])^{(4)}
\]
\[
\oplus (K[W_2(2)])^{(6)}\otimes (K[W_2(3)])^{(6)}
\oplus (K[W_2(2)])^{(4)}\otimes (K[W_2(2^2)])^{(8)}
\]
\[
\oplus (K[W_2(2)])^{(2)}\otimes (K[W_2(3)])^{(6)}\otimes (K[W_2(2^2)])^{(4)}
\]
\[
\oplus (K[W_2(3)])^{(12)}
\oplus (K[W_2(2^2)])^{(12)}
\]
and obtain the multiplicity $m(6^2)$.
We shall mention that $W_2(6^2)$ participates with multiplicity one
in the components
\[
(K[W_2(2)])^{(12)}, (K[W_2(2)])^{(8)}\otimes (K[W_2(2^2)])^{(4)},
(K[W_2(2)])^{(4)}\otimes (K[W_2(2^2)])^{(8)},
\]
\[
(K[W_2(2)])^{(2)}\otimes (K[W_2(3)])^{(6)}\otimes (K[W_2(2^2)])^{(4)},
(K[W_2(3)])^{(12)}, (K[W_2(2^2)])^{(12)}
\]
and with multiplicity 2 in the component
$(K[W_2(2)])^{(6)}\otimes (K[W_2(3)])^{(6)}$.
\end{proof}

\begin{lemma}\label{generators for 66}
The elements $w_1,w_2,w_3',w_3'',w_4,w_5,w_6,w_7$
given in {\rm (\ref{w1 and the others}), (\ref{w5}),
(\ref{w6}), (\ref{w31})} and  {\rm (\ref{w32})} are highest weight vectors
of submodules $W_2(6^2)$ and are linearly independent in $S$.
\end{lemma}

\begin{proof}
It is a direct checking to see that the eight elements are linearly independent.
All these elements are of degree (6,6).
In order to see that they are highest weight vectors it is sufficient to
verify that each $w_i(x,y)$ has the property $w_i(x\vert y,x)=0$, as in the proof of
Lemma \ref{lemma-B}.
\end{proof}

We think that, for possible applications,
it is important to know how we have found the elements
$w_1,w_2$, $w_3',w_3''$, $w_4,w_5,w_6,w_7$.
The element $v$ from (\ref{v}) belongs to $W_2(2^2)$ and is a highest
weight vector. Similarly, $u$ from (\ref{u}) is a highest weight
vector in $W_2(2^2)\subset (K[W_2(2)])^{(4)}$. This explains why
$w_1,w_2,w_4$ and $w_7$ are highest weight vectors in the submodules
$W_2(6^2)$ of $(K[W_2(2)])^{(12)}$, $(K[W_2(2)])^{(8)}\otimes (K[W_2(2^2)])^{(4)}$,
$(K[W_2(2)])^{(4)}\otimes (K[W_2(2^2)])^{(8)}$ and
$(K[W_2(2^2)])^{(12)}$, respectively.
The tensor product
$(K[W_2(2)])^{(2)}\otimes (K[W_2(3)])^{(6)}\otimes (K[W_2(2^2)])^{(4)}$
contains one submodule $W_2(6^2)$. Since $(K[W_2(2^2)])^{(4)}\cong W_2(2^2)$
is one-dimensional and spanned by the highest weight vector $v$, we need to find a
highest weight vector of the submodule $W_2(4^2)$ of
$(K[W_2(2)])^{(2)}\otimes (K[W_2(3)])^{(6)}$. It has been done by the method
of Lemma \ref{lemma-B}. We are looking for a linear combination $w_{(4^2)}(x,y)$
of products $\text{\rm tr}(z_1z_2)\text{\rm tr}(z_3z_4z_5)\text{\rm tr}(z_6z_7z_8)$,
$z_i=x,y$, which is homogeneous of degree $(4,4)$ and satisfies the condition
$w_{(4^2)}(x\vert y,x)=0$. All possible summands are
\[
\text{\rm tr}(x^2)\text{\rm tr}(x^2y)\text{\rm tr}(y^3),
\text{\rm tr}(x^2)\text{\rm tr}^2(xy^2),
\]
\[
\text{\rm tr}(xy)\text{\rm tr}(x^3)\text{\rm tr}(y^3),
\text{\rm tr}(xy)\text{\rm tr}(x^2y)\text{\rm tr}(xy^2),
\]
\[
\text{\rm tr}(y^2)\text{\rm tr}(x^3)\text{\rm tr}(xy^2),
\text{\rm tr}(y^2)\text{\rm tr}^2(x^2y).
\]
We determine the coefficients of these summands from the condition
$w_{(4^2)}(x\vert y,x)=0$. It has turned out that $w_{(4^2)}$ has a determinantal form
and in this way we obtain $w_5$. The highest weight vector
$w_6\in (K[W_2(3)])^{(12)}$ is obtained similarly as a linear combination
of homogeneous products of four $\text{\rm tr}(z_{i_1}z_{i_2}z_{i_3})$,
of degree (6,6), $z_i=x,y$.
Finally, the proof of Lemma \ref{multiplicities for 33 and 66} gives that
the multiplicity of $W_2(6^2)$ in $(K[W_2(2)])^{(6)}\otimes (K[W_2(3)])^{(6)}$
is equal to 2. Since
\[
(K[W_2(2)])^{(6)}\cong (K[W_2(3)])^{(6)} \cong W_2(6)\oplus W_2(4,2),
\]
the Littlewood-Richardson rule gives that we may fix one of the highest
weight vectors of the two copies of $W_2(6^2)$, namely $w_3'$, in
$W_2(4,2)\otimes W_2(4,2)$, and the second, namely $w_3''$, in
$W_2(6)\otimes W_2(6)$. For $w_3'$ we use the same arguments as for
$w_5$, replacing $v$ with $u\in (K[W_2(2)])^{(4)}$. For $w_3''$, from
the expression (\ref{new definition of w32})
and the equalities
$\delta^7(\text{\rm tr}^3(y^2))=
\delta^7(\text{\rm tr}^2(y^3))=0$ we obtain
$w_3''(x\vert y,x)=\delta(w_3(x,y))=0$.

\begin{lemma}\label{search expression for ww}
There exist constants $\xi_1,\xi_2,\xi_3',\xi_3'',\xi_4,\xi_5,\xi_6,\xi_7$
such that the relation
\begin{equation}\label{relation with unknown coefficients}
w^2-(\xi_1w_1+\xi_2w_2+\xi_3'w_3'+\xi_3''w_3''+\xi_4w_4+\xi_5w_5+\xi_6w_6+\xi_7w_7)=0
\end{equation}
holds in $C_{32}$.
\end{lemma}

\begin{proof}
By Lemma \ref{the modules generated by the coefficients},
the coefficients $a$ and $b$ in (\ref{the relation for w})
belong to $GL_2$-submodules of the algebra $S$
isomorphic to $W_2(3^2)$ and $W_2(6^2)$, respectively.
Lemma \ref{multiplicities for 33 and 66} gives that $S$ does not contain
a submodule $W_2(3^2)$. Hence the coefficient $a$ is equal to 0 and
$w$ satisfies a relation $w^2+b=0$ for some $b\in W_2(6^2)\subset S$.

Representation theory of $GL_d$ gives that if $W_i$, $i=1,\ldots,k$,
are $k$ isomorphic copies of $W_d(\lambda)$ and $w_i\in W_i$ are highest weight
vectors, then the highest weight vector of any submodule $W_d(\lambda)$
of the direct sum $W_1\oplus\cdots\oplus W_k$ has the form
$\xi_1w_1+\cdots+\xi_kw_k$ for some $\xi_i\in K$.
Any $k$ linearly independent highest weight vectors can serve
as a basis of the $GL_d$-module $W_1\oplus\cdots\oplus W_k$.
By Lemma \ref{multiplicities for 33 and 66}, $W_2(6^2)$ participates in
$S$ with multiplicity 8 and
Lemma \ref{generators for 66} gives such 8 linearly
independent highest weight vectors. Hence $w$ satisfies the desired relation.
\end{proof}

Now the problem is to find the coefficients
$\xi_1,\xi_2,\xi_3',\xi_3'',\xi_4,\xi_5,\xi_6,\xi_7$ in the relation
(\ref{relation with unknown coefficients}). We have proceeded in the following way,
using Maple for the concrete calculations.

We have evaluated (\ref{relation with unknown coefficients}) for
various traceless matrices. First, we assume that
$x$ is as in (\ref{x is diagonal}) and
\[
y=\left( \begin{array}{ccc}
                0 & 1 & 0 \\
                0 & 0 & 1 \\
                1 & 0 & 0 \\
             \end{array} \right).
\]
The only nonzero expressions are for $w, w_3''$ and $w_6$.
Comparing the coefficients of $x_1^6$ and $x_1^4x_2^2$, we obtain that
\[
4-360\xi_3''=0\quad \text{\rm and}\quad -3-(2160 \xi_3''+81 \xi_6)=0,
\]
which gives that
\begin{equation}\label{first equation on xi}
\xi_3''=1/90,\quad \xi_6=-1/3.
\end{equation}
The second step is to use
\[
x=\left( \begin{array}{ccc}
                1 & 0 & 0 \\
                0 & -1 & 0 \\
                0 & 0 & 0 \\
             \end{array} \right),\quad
y=\left( \begin{array}{ccc}
                0 & 0 & 0 \\
                0 & 1 & 0 \\
                0 & 0 & -1 \\
             \end{array} \right),
\]
which gives $w=0$ and $27\xi_1-9\xi_3'+36\xi_3''-3\xi_6=0$. Hence, by
(\ref{first equation on xi}),
\begin{equation}\label{second equation on xi}
  27\xi_1 - 9 \xi_3' = -7/5.
\end{equation}
Then, again we fix $x$ as in (\ref{x is diagonal}) and
\[
y=\left( \begin{array}{ccc}
                0 & 1 & 0 \\
                1 & 0 & 0 \\
                0 & 0 & 0 \\
             \end{array} \right).
\]
We obtain $w=0$ and
\[
(64\xi_1+16\xi_2-8\xi_3'+72\xi_3''+4\xi_4-2\xi_5+\xi_7)x_1^6
\]
\[
+(192\xi_1-56\xi_3'+144\xi_3''-12\xi_4-8\xi_5-12\xi_6-6\xi_7)x_1^5x_2
\]
\[
+(384\xi_1-136\xi_3'+504\xi_3''+12\xi_4+2\xi_5-48\xi_6+15\xi_7)x_1^4x_2^2
\]
\[
+(448\xi_1-32\xi_2-176\xi_3'+864\xi_3''-8\xi_4+16\xi_5-72\xi_6-20\xi_7)x_1^3x_2^3
\]
\[
+(384\xi_1-136\xi_3'+504\xi_3''+12\xi_4+2\xi_5-48\xi_6+15\xi_7)x_1^2x_2^4
\]
\[
+(192\xi_1-56\xi_3'+144\xi_3''-12\xi_4-8\xi_5-12\xi_6-6\xi_7)x_1x_2^5
\]
\[
+(64\xi_1+16\xi_2-8\xi_3'+72\xi_3''+4\xi_4-2\xi_5+\xi_7)x_2^6=0.
\]
Since the coefficients of $x_1^ix_2^{6-i}$ are equal to 0, this gives four more relations
for the $\xi$'s. Combining with (\ref{first equation on xi}) and
(\ref{second equation on xi}), we have
\begin{equation}\label{third group equations on xi}
\begin{array}{c}
\xi_1= -1/54+1/2\xi_4+3/4\xi_7,\quad
\xi_2= 1/54-3/2\xi_4-7/4\xi_7, \\
\\
\xi_3'= 1/10+3/2\xi_4+9/4\xi_7,\quad
\xi_5 = - 4/9 + 3/2 \xi_7. \\
\end{array}
\end{equation}
Finally, using $x$ from (\ref{x is diagonal}) and $y$ from (\ref{matrices}),
together with the expressions for the $\xi$'s from
(\ref{first equation on xi}), (\ref{second equation on xi}) and
(\ref{third group equations on xi}),
we annulate the coefficient of
$x_1^3x_2^3y_{11}^3y_{12}y_{23}y_{31}$ to derive
\[
-27\xi_4-81/2 \xi_7  + 3=0.
\]
Similarly, the coefficient of $x_1x_2^5y_{12}y_{22}y_{23}^2y_{31}y_{32}$ gives
\[
36\xi_4+90\xi_7+4/3=0.
\]
Hence the only possibility is
\[
\xi_4=1/3.\quad \xi_7=-4/27
\]
and we conclude that the relation is the given in
(\ref{our defining relation}).

\section*{Acknowledgements}

This project was started when the second author
visited the Department of Mathematics of the
National University of Singapore. He is very grateful for the hospitality
and the creative atmosphere during his stay in Singapore.

\end{document}